\documentclass[11pt]{amsart}
\usepackage{amstext,amssymb,amsmath,amsbsy,dsfont,tikz}

 \usepackage[left=2.65cm,right=2.65cm,top=3.2cm,bottom=3.2cm]{geometry}

\usepackage{amstext,amssymb,amsmath,amsbsy}

\usepackage{tikz}
\usepackage{hyperref,cleveref}
\usepackage{amscd}
\usepackage{amsfonts}
\usepackage{indentfirst}   
\usepackage{verbatim}
\usepackage{amsmath}
\usepackage{amsthm}
\usepackage{enumerate}
\usepackage{graphicx}
\usepackage{subfig}
\usepackage{color}
\usepackage[utf8]{inputenc}
\usepackage[T1]{fontenc}     % un second package
\usepackage[english]{babel} 
\usepackage{amssymb}
\usepackage{bm}
\usepackage{ifthen, soulutf8}

\usepackage{mathtools}
\usepackage{etoolbox}

\DeclareFontFamily{U}{matha}{\hyphenchar\font45}
\DeclareFontShape{U}{matha}{m}{n}{
<-6> matha5 <6-7> matha6 <7-8> matha7
<8-9> matha8 <9-10> matha9
<10-12> matha10 <12-> matha12
}{}
\DeclareSymbolFont{matha}{U}{matha}{m}{n}

\DeclareFontFamily{U}{mathx}{\hyphenchar\font45}
\DeclareFontShape{U}{mathx}{m}{n}{
<-6> mathx5 <6-7> mathx6 <7-8> mathx7
<8-9> mathx8 <9-10> mathx9
<10-12> mathx10 <12-> mathx12
}{}
\DeclareSymbolFont{mathx}{U}{mathx}{m}{n}

\DeclareMathDelimiter{\vvvert} {0}{matha}{"7E}{mathx}{"17}%

\DeclarePairedDelimiterX{\normiii}[1]
{\vvvert}
{\vvvert}
{\ifblank{#1}{\:\cdot\:}{#1}}

\newtheorem{theorem}{Theorem}

\newtheorem{lemma}{Lemma}

\newtheorem{proposition}{Proposition}

\newtheorem{remark}{Remark}

\newtheorem{definition}{Definition}
      \newcommand{\mC}{\mathbb{C}}

\numberwithin{equation}{section}
\numberwithin{theorem}{section}
\numberwithin{prop}{section}
\numberwithin{lem}{section}
\numberwithin{definition}{section}
\numberwithin{rem}{section}
\numberwithin{lemma}{section}
  \numberwithin{proposition}{section}

\newcommand{\bH}{{\bf H}}
\newcommand{\cE}{{\mathcal E}}
\newcommand{\cH}{{\mathcal H}}

\newcommand{\tlambda}{\widetilde \lambda}

\newcommand{\hE}{\hat E}
\newcommand{\hH}{\hat H}

\newcommand{\hcE}{\hat {\mathcal  E}}
\newcommand{\hcH}{\hat {\mathcal H}}

\newcommand{\cN}{{\mathcal N}}

      \newcommand{\hmu}{\hat \mu}
      \newcommand{\heps}{\hat \eps}
      \newcommand{\hJ}{\hat J}

\newcommand{\cT}{\mathcal T}

\newcommand{\dive}{\operatorname{div}}
\newcommand{\curl}{\operatorname{curl}}

\newcommand{\eps}{\varepsilon}
\renewcommand{\epsilon}{\varepsilon}

\newcommand{\mN}{\mathbb{N}}

\newcommand{\mR}{\mathbb{R}}
\newcommand{\R}{\mathbb{R}}

\newcommand{\C}{\mathbb{C}}

\newcommand{\hu}{\hat u}

\newcommand{\hv}{\hat v}

\newcommand{\M}{{\mathcal M}}
\newcommand{\hxi}{\hat \xi}

\newcommand{\cJ}{{\mathcal J}}

\newcommand{\bcJ}{{\bf {\mathcal J}}}

\newcommand{\bJ}{{\bf J}}

\newcommand{\hbJ}{\hat {\bf J}}

\newcommand{\hcJ}{\hat {\mathcal J}}

\newboolean{showcomments}
\setboolean{showcomments}{true}
\ifthenelse{\boolean{showcomments}}
{ \newcommand{\mynote}[3]{
    \fbox{\bfseries\sffamily\scriptsize#1}
    {\small$\blacktriangleright$\textsf{\textit{\color{#3}{#2}}}$\blacktriangleleft$}}}
{ \newcommand{\mynote}[3]{}}
\setul{0.5ex}{0.3ex}
\setulcolor{red}

\title[The transmission eigenvalue problem for Maxwell equations]{The completeness of the generalized eigenfunctions
and an upper bound for the counting function of the transmission eigenvalue problem for Maxwell equations}

\author{Jean Fornerod}
\address[Jean Fornerod]{Ecole Polytechnique F\'ed\'erale de Lausanne, EPFL,  SB, CAMA,
	\newline\indent Station 8,  CH-1015 Lausanne, Switzerland.}
\email{jean.fornerod@epf.ch}

\author{Hoai-Minh Nguyen}
\address[Hoai-Minh Nguyen]{Ecole Polytechnique F\'ed\'erale de Lausanne, EPFL,  SB, CAMA, 
\newline\indent Station 8,  CH-1015 Lausanne, Switzerland.}
\email{hoai-minh.nguyen@epfl.ch}

\begin{document}

\maketitle 

\begin{abstract}  Cakoni and Nguyen recently proposed very general conditions on the coefficients of Maxwell equations for which they established the discreteness of the set  of eigenvalues of the transmission problem and studied their locations. In this paper, we establish  the completeness of the generalized eigenfunctions and derive an optimal upper bound for the counting function under these conditions, assuming additionally that the coefficients are twice continuously differentiable. The approach is based on the spectral theory of Hilbert-Schmidt operators. 

\end{abstract}

\section{Introduction}

Let $\Omega \subset \R^3$ be a bounded domain of class $C^3$. Let $\epsilon, \, \mu, \, \hmu, \, \heps \in [L^\infty(\Omega)]^{3 \times 3}$ be symmetric and uniformly elliptic.  
A complex number $\omega \in \C$ is called a \textit{transmission eigenvalue} if there exists a non-zero solution $(E,H,\hat{E},\hat{H}) \in [L^2(\Omega)]^{12}$ of the following Cauchy problem 
\begin{equation}
\label{e11}
\left \{
\begin{array}{cl}
\nabla \times E  =  i \omega \mu  H & \mbox{ in } \Omega,  \\[6pt]
\nabla \times H  = - i \omega \eps E  & \mbox{ in } \Omega,\\
\end{array}
\right .
\quad \left \{
\begin{array}{cl}
\nabla \times \hat{E}  = i \omega \hat{\mu} \hat{H} & \mbox{ in } \Omega, \\[6pt]
\nabla \times \hat{H}  = - i \omega \heps \hat{E} & \mbox{ in } \Omega,  
\end{array}
\right . 
\end{equation}
\begin{equation}
\label{e12}
(\hat{E}-E)\times \nu = 0  \mbox{ on } \partial \Omega,  \quad \text{ and } \quad  (\hat{H}-H)\times \nu = 0 \: \text{ on }\partial \Omega. 
\end{equation}
Here and in what follows,  $\nu$ denotes the unit, outward, normal vector to $\partial \Omega$. 

The  transmission eigenvalue problem, proposed by Kirsch \cite{KIrsch86} and Colton and Monk \cite{CM88},  has been an active research topic  in the   inverse scattering theory for inhomogeneous media. It has  a  connection with the injectivity  of the relative scattering operator. Transmission eigenvalues are  related to interrogating frequencies for which there is an incident field that does not scatterer by the medium. We refer the reader to \cite{CCH16} for a recent and self-contained introduction to the topic.

Cakoni and Nguyen \cite{Cakoni-Ng21} have recently studied 
the transmission problem for Maxwell equations in a very general setting. Under the assumption $\eps, \, \mu, \, \heps, \, \hat{\mu}$ are of class $C^1$ in a neighborhood of  the boundary,  they proposed  the following  condition: 
\begin{equation}\label{H}
 \mbox{$\eps, \, \mu,  \, \heps, \, \hat{\mu}$ are isotropic  on $\partial \Omega$, and 
$\epsilon \neq \hat{\epsilon}, \: \: \mu \neq \hat{\mu}, \: \: \epsilon/\mu \neq  \hat{\epsilon}/ \hat{\mu} \: \: \text{ on }\partial \Omega$}  
\end{equation}
(see \Cref{rem-iso} for the convention used in \eqref{H}). Under this assumption,  Cakoni and Nguyen  showed that the set of eigenvalues $\lambda_j$ of system \eqref{e11}-\eqref{e12} is discrete.  They also studied the location of the eigenvalues under this condition. More precisely, they showed that, for every $\gamma >0$, there exists $\omega_0 > 0$ such that if $\omega \in \mC$ with $|\Im(\omega^2)| \ge \gamma |\omega|^2$ and $|\omega| \ge \omega_0$, then $\omega$ is not a transmission eigenvalue. Their analysis is inspired and guided by the famous work of  Agmon, Douglis, and Nirenberg \cite{ADNI, ADNII} on complementing  boundary conditions. 

In this paper, we study further spectral properties of the transmission problem under assumption \eqref{H}  given above. More precisely, we establish  the completeness  of the generalized eigenfunctions and  derive an optimal upper bound for the counting function of the transmission eigenvalues. 

Before stating our results, as in \cite{Cakoni-Ng21},  we denote
\begin{multline}
\bH (\Omega) : =\Big\{ (u,v,\hat{u},\hat{v}) \in [L^2(\Omega)]^{12} : \dive (\eps u)= \dive (\mu v) = \dive (\heps \hat{u}) = \dive (\hat{\mu} \hat{v}) = 0 \text{ in }\Omega,  \\[6pt] \mbox{ and } \quad 
\heps \hu \cdot \nu - \eps u \cdot \nu = \hmu \hv \cdot \nu - \mu v  \cdot \nu = 0 \text{ on }\partial \Omega \Big\}. 
\end{multline}
The functional space $\bH(\Omega)$, which plays a role in the analysis in \cite{Cakoni-Ng21} as well as in this paper,  is  a Hilbert space with the standard $[L^2(\Omega)]^{12}$-scalar product.  One of the motivations for the definition of $\bH(\Omega)$ is  the fact that if $(E, H, \hE, \hH) \in [L^2(\Omega)]^{12}$ is an eigenfunction of the transmission problem, i.e., a  solution of \eqref{e11} and \eqref{e12} for some $\omega \in \mC$, then $(E, H, \hE, \hH) \in \bH(\Omega)$ except for $\omega = 0$. The other motivation is on the compactness of $\cT_k$ defined below. 

\medskip 
The first main result of this paper is  on the completeness of the generalized eigenfunctions. We have 

\begin{theorem}
\label{thm1}
Assume that  $\eps, \, \mu, \, \heps, \, \hat{\mu}  \in [C^2(\bar \Omega)]^{3 \times 3}$ and \eqref{H} holds.   The space spanned by the generalized eigenfunctions  is complete in $\bm{H}(\Omega)$, i.e., the space spanned by them is dense in  $\bH(\Omega)$. 
\end{theorem}

\begin{remark} \rm See also \Cref{rem-Form} for a discussion on another version of \Cref{thm1}. 
\end{remark}

The second main result of this paper is  on  an upper  bound for the counting function $\cN$. This function is defined by, for $t > 0$, 
\begin{equation}\label{def-N}
\cN(t) := \# \Big\{ j : |\lambda_j| \leq t \Big\}.
\end{equation}
Concerning the behavior of $\cN(t)$ for large $t$, we have 

\begin{theorem} \label{thm2} Assume that  $\eps, \, \mu, \, \heps, \, \hat{\mu}  \in [C^2(\bar \Omega)]^{3 \times 3}$ and \eqref{H} holds. There exists a constant $c>0$  such that,  for $t > 1$, 
\begin{equation}
\cN(t) \leq ct^3.
\end{equation}
\end{theorem}

As a consequence of \Cref{thm1}, the set of the transmission eigenvalues is infinite. This fact is new as far as we know. 
\Cref{thm2}, complement to \Cref{thm1},  gives an upper bound for the density of the distribution of the transmission eigenvalues. 
This upper bound is optimal 
in the sense that it has the same order as the standard Weyl laws for the Maxwell equations \cite{Weyl12,BF07}. 

\medskip 
Some comments on \Cref{thm1} and \Cref{thm2}  are in order. The generalized eigenfunctions associated with $\lambda_j$, considered in \Cref{thm1},  are understood as the generalized eigenfunctions of the operator $\cT_k$,  defined in \eqref{def-Tk} below, corresponding to the eigenvalue $ (i \lambda_j  - k)^{-1} $ of $\cT_k$. One can show that it is independent of $k$ as long as $\cT_k$ is well-defined (and compact).  In the conclusion of  \Cref{thm2}, the multiplicity of eigenvalues is taken into account. The meaning of the multiplicity $\lambda_j$ is understood as the multiplicity of the eigenvalue  $ (i \lambda_j  - k)^{-1}$ of the operator $\cT_k$. Again, this is independent of $k$.  These facts follow from \cite[Theorem 12.4]{Agmon} after using \Cref{lem-Tmod} on the modified resolvent of $\cT_k$. The multiplicity and the generalized eigenfunctions corresponding to $\lambda_j$ are then understood as the multiplicity of $(i \lambda_j - k)^{-1}$ and the generalized eigenfunctions corresponding to $(i \lambda_j - k)^{-1}$ both corresponding to $\cT_k$ from now on.

\medskip
We recall here the definition of a generalized eigenfunction and the multiplicity of its corresponding eigenvalue, see e.g. \cite[Definition 12.5]{Agmon}, for the convenience of the reader. 

\begin{definition} \label{def-GE}
Let $A: H \to H$ be a linear and bounded operator on a Hilbert space $H$. Let $\lambda$ be an eigenvalue of $A$. An element $v \in H \setminus \{0 \}$ is a called a {\rm generalized eigenfunction} of $T$ if there exists a positive integer $m$ such that
\begin{equation}
(\lambda-A)^m v = 0.
\end{equation}
The  multiplicity  of the eigenvalue $\lambda$  is defined as the dimension of the set $\bigcup_{m \in \mN_*} {\rm Ker}(\lambda - A)^m. $
\end{definition}

The study of the  transmission eigenvalue problem for Maxwell's equations is not as complete as  for the scalar case, which is discussed briefly below. Before \cite{Cakoni-Ng21}, the  discreteness results can be found in \cite{Haddar04,CHM15} (see also \cite{Chesnel12}) where the case of $\mu = \heps = \hmu = I $, and $\eps - I $ invertible in a neighborhood  $\partial \Omega$ was considered. Concerning the other aspects, Cakoni, Gintides, and Haddar \cite{CGH10} studied the existence of real transmission eigenvalues,  and  Haddar and Meng \cite{HM18} studied the completeness of eigenfunctions for the setting related to the one in \cite{CHM15} mentioned above.  In the isotropic case, under the assumption  $\mu = \hmu$ and $\eps \mu \neq \heps \hmu$, Vodev recently derived  a parabolic eigenvalue-free region  \cite{Vodev21}.

The structure of spectrum of the transmission eigenvalue problem is better understood in the case of scalar inhomogeneous Helmoltz equations. In this case, the  transmission eigenvalue problem can be stated as follows. Let $d \ge 2$ and $\Omega$ be an open,  bounded  Lipschitz subset of $\mR^d$. Let $A_1$ and $A_2$ be two ($d \times d$) symmetric, uniformly elliptic,   matrix-valued  functions and $\Sigma_1$ and $\Sigma_2$ be two bounded positive functions both defined on $\Omega$.  A complex number $\omega$ is called an eigenvalue of the  transmission eigenvalue problem, or a transmission eigenvalue, if there exists a non-zero solution $(u_1, u_2) \in [H^1(\Omega)]^2$ of the system 
 \begin{equation}\label{pro1a}  
 \left\{\begin{array}{lll}
 \dive(A_1 \nabla u_1) + \omega^2 \Sigma_1 u_1= 0 ~~&\text{ in}~\Omega, \\[6pt]
 \dive(A_2 \nabla u_2) + \omega^2 \Sigma_2 u_2= 0 ~~&\text{ in}~\Omega, 
 \end{array} \right. 
 \end{equation}
 \begin{equation}\label{pro1b}
  u_1 =u_2, \quad  A_1 \nabla u_1\cdot \nu = A_2 \nabla u_2\cdot \nu ~\text{ on } \partial \Omega.  
 \end{equation}
The discreteness of transmission eigenvalues for  the Helmholtz equation has been  investigated extensively in the literature. The state-of-the-art results  on the discreteness of transmission eigenvalues are given in \cite{MinhHung17}. More specifically, it was shown in \cite{MinhHung17} that the transmission eigenvalue problem has a discrete spectrum if the coefficients are smooth  only near the boundary, and 
\begin{enumerate}
\item[i)] $A_1(x),  \, A_2(x)$ satisfy the complementing boundary condition with respect to  $\nu(x)$ for all $x \in \partial \Omega$, i.e., 
for all $x \in \partial \Omega$ and for all $\xi \in \mR^d \setminus \{0\}$ with  $\xi \cdot \nu = 0$, we have 
\begin{equation*}
(A_2 \nu \cdot \nu ) ( A_2 \xi \cdot \xi )  - ( A_2 \nu \cdot  \xi )^2 \neq  ( A_1 \nu, \nu )(  A_1 \xi \cdot  \xi )  - ( A_1 \nu \cdot \xi )^2, 
\end{equation*}
\item[ii)]  $( A_1 \nu \cdot \nu ) \Sigma_1 \neq (  A_2 \nu \cdot \nu ) \Sigma_2$ for all $x \in \partial \Omega$. 
\end{enumerate}
Additional results in \cite{MinhHung17} also include various combinations of the sign of contrasts $A_1-A_2$ and $\Sigma_1-\Sigma_2$ on the boundary. Assume i) and ii) and $A_1, A_2, \Sigma_1, \Sigma_2$ are continuous in $\bar \Omega$,  the Weyl laws for eigenvalues and the completeness of the generalized eigenfunctions in $[L^2(\Omega)]^2$ were recently established by Nguyen and (Q. H.) Nguyen \cite{MinhHung2}.    Previous results on discreteness can be found in  \cite{BCH11,  LV12, Sylvester12} and references therein.  
Completeness of  transmission eigenfunctions and estimates on the  counting function were  studied  by Robbiano  \cite{Robbiano13,Robbiano16}  for $C^\infty$  boundary and coefficients,  and for the case $A_1 = A_2 = I$. Again in $C^\infty$ isotropic setting, Vodev \cite{Vodev15}, \cite{Vodev18} proved the sharpest known results  on eigenvalue free zones and Weyl's law with an estimate for the remainder.

The Cauchy problem also naturally appears in the context of negative index materials after using reflections as initiated in \cite{Ng-Complementary} (see also \cite{Ng-Superlensing-Maxwell}).  The well-posedness and the limiting absorption principle for the Helmholtz equation with sign-changing coefficients
were developed by Nguyen \cite{Ng-WP} using the Fourier and multiplier approach. Similar problems for the Maxwell equations were studied by Nguyen and Sil \cite{NgSil}. Both papers \cite{Ng-WP},  \cite{NgSil}  deal with the stability question of negative index materials,  and are the starting point for the analysis of the  transmission eigenvalue problems  in \cite{MinhHung17, MinhHung2, Cakoni-Ng21}.  Other aspects and applications of negative index materials as well as the stability and instability the Cauchy problem \eqref{e11} and \eqref{e12} are discussed in \cite{Ng-Superlensing-Maxwell,  Ng-Negative-Cloaking-M, Ng-CALR-M, Ng-CALR-O-M} and the references therein.

\medskip 
The key and the starting point of the analysis in \cite{Cakoni-Ng21} is the following result \cite[Propositions 4.1 and 4.2]{Cakoni-Ng21}:
\begin{theorem}[Cakoni \& Nguyen] \label{thm-CN} 
Assume that  $\eps, \, \mu, \, \heps, \, \hat{\mu}  \in [C^1(\bar \Omega)]^{3 \times 3}$ and \eqref{H} holds,  and  let $\gamma > 0$. There exist two  constants $k_0 \ge 1$ and $C > 0$ such that for $k \in \mC$ with $|\Im{(k^2)}| \ge \gamma |k|^2$ and $|k| \ge k_0$,  for every
 $(J_e, J_m, \hJ_e, \hJ_m) \in [L^2(\Omega)]^{12}$, there exists   a unique solution   $(E, H, \hE, \hH) \in [L^2(\Omega)]^{12}$ of 
\begin{equation}\label{sys-ITE-1}
\left\{\begin{array}{c}
\nabla \times E = k \mu H + J_e \mbox{ in } \Omega, \\[6pt]
\nabla \times H = - k \eps E + J_m \mbox{ in } \Omega, 
\end{array}\right. \quad \left\{\begin{array}{c}
\nabla \times \hE = k \hmu \hH + \hJ_e \mbox{ in } \Omega, \\[6pt]
\nabla \times \hH = -  k \heps \hE + \hJ_m \mbox{ in } \Omega,  
\end{array}\right.
\end{equation}
\begin{equation}\label{bdry-ITE-1}
(\hE - E) \times \nu = 0 \mbox{ on } \partial \Omega, \quad \mbox{ and } \quad (\hH - H) \times \nu = 0 \mbox{ on } \partial \Omega. 
\end{equation}
Moreover, if  $(J_e, J_m, \hJ_e, \hJ_m) \in [H(\dive, \Omega)]^4$ with $(J_{e}\cdot \nu - \hJ_{e} \cdot \nu, J_{m} \cdot \nu -  \hJ_{m} \cdot \nu) \in [H^{1/2}(\partial \Omega)]^2$, then 
\begin{multline}\label{pro-main-1-est}
 |k| \,  \| (E, H, \hE, \hH) \|_{L^2(\Omega)}  + \| (E, H, \hE, \hH) \|_{H^1(\Omega)} 
 \le C \| (J_e, J_m, \hJ_e, \hJ_m)\|_{L^2(\Omega)} \\[6pt]
 + \frac{C}{|k|}  \| (\dive J_e,\dive  J_m, \dive \hJ_e,\dive  \hJ_m)\|_{L^2(\Omega)} +  \frac{C}{|k|}  \| (J_{e} \cdot \nu - \hJ_e \cdot \nu, J_m \cdot \nu -  \hJ_m \cdot \nu)\|_{H^{1/2}(\partial \Omega)}  . 
\end{multline}
\end{theorem}

\begin{remark}  \rm In \cite{Cakoni-Ng21}, the coefficients are only assumed to be of class $C^1$ near the boundary and  a variant of \eqref{pro-main-1-est},  where the $\| \cdot \|_{H^1(\Omega)}$ is replaced by $\| \cdot \|_{H^1(D \cap \Omega)}$ for some neighborhood $D$ of $\partial \Omega$ (see \cite[(4.4) of Proposition 4.1]{Cakoni-Ng21}), was established. Nevertheless, under the smoothness assumption considered here, \eqref{pro-main-1-est} follows immediately by the same analysis. 
\end{remark}

Fix  $k \in \mC$ such that the conclusions in \Cref{thm-CN} hold.  One can then define the operator $\cT_k$ as follows 
\begin{equation}\label{def-Tk}
\begin{array}{rclc}
\cT_k: & \bm{H}(\Omega) & \to & \bm{H}(\Omega) \\[6pt]
& (\cJ_e, \cJ_m, \hcJ_e, \hcJ_m) & \mapsto & (E,H,\hat{E},\hat{H}), 
\end{array}
\end{equation}
where $(E, H, \hE, \hH)$ is the unique solution of, with $(J_e, J_m, \hJ_e, \hJ_m) = (\mu \cJ_m, - \eps \cJ_e, \hmu \hcJ_m, -  \heps \hcJ_e)$, 
\begin{equation}
\label{e11p}
\left \{
\begin{array}{cl}
\nabla \times E  = k  \mu  H  + J_e & \mbox{ in } \Omega,  \\[6pt] 
\nabla \times H  = - k \eps E + J_m  & \mbox{ in } \Omega,
\end{array}
\right .
\quad \left \{
\begin{array}{cl}
\nabla \times \hE  = k  \hmu  \hH  + \hJ_e & \mbox{ in } \Omega,  \\[6pt] 
\nabla \times \hH  = - k \heps \hE +  \hJ_m  & \mbox{ in } \Omega,
\end{array}
\right . 
\end{equation}
\begin{equation}
\label{e12p}
(\hat{E}-E)\times \nu = 0  \mbox{ on } \partial \Omega,  \quad \text{ and } \quad  (\hat{H}-H)\times \nu = 0 \: \text{ on }\partial \Omega. 
\end{equation}
From \eqref{pro-main-1-est} and the compactness criterion related to the Maxwell equations, one can derive that  $\cT_k$ is compact. It is easy to check that $\omega$ is an eigenvalue of the transmission problem if and only if $(i \omega - k)^{-1}$ is an eigenvalue of $\cT_k$. The discreteness of the eigenvalues of the transmission problem then follows from the discreteness of the eigenvalues of $\cT_k$.

In this paper,  to derive further spectral properties of the transmission problem, we develop the analysis in \cite{Cakoni-Ng21} in order to be able to  apply the spectral theory of Hilbert-Schmidt operators. This strategy was previously used in the acoustic setting \cite{MinhHung2}.  To this end, we establish a regularity result (see \Cref{thm-reg}) for solutions given in \Cref{thm-CN}. In addition to this, one of the main ingredients in the proof of  \Cref{thm1} is  the density of the range of the map $\cT_k$ in $\bH(\Omega)$ with respect to the $[L^2(\Omega)]^{12}$-norm (see \Cref{pro-density}).  The  proof of \Cref{thm1} is also given in a way which does not involve any extra topological property of $\Omega$ than its connectivity (see Step 2 of the proof of \Cref{pro-density}). 

\medskip 
The paper is organized as follows. In \Cref{sect-reg}, we establish the regularity result on the transmission problem.   The last two sections are devoted to the proof of \Cref{thm1} and \Cref{thm2}, respectively.

\section{A regularity result for the transmission problem}\label{sect-reg}

The following regularity result for the Maxwell transmission problem is the main result of this section (compare with \Cref{thm-CN}).

\begin{theorem} \label{thm-reg}
Let  $\eps, \, \mu, \, \heps, \, \hat{\mu} \in [C^2(\bar{\Omega})]^{3 \times 3}$ be symmetric,  and let $\gamma > 0$.  Assume that there exist $\Lambda\geq 1$ and  $\Lambda_1>0$ such that 
\begin{equation}\label{thm-reg-as1}
\Lambda^{-1} \leq \epsilon,\mu,\hat{\epsilon},\hat{\mu} \leq \Lambda  \text{ in }\Omega,  \quad   \|(\epsilon,\mu,\hat{\epsilon},\hat{\mu})\|_{C^2(\bar{\Omega})} \leq \Lambda, 
\end{equation}
\begin{equation}\label{thm-reg-as2}
\eps, \, \mu, \, \heps, \, \hat{\mu} \text{ are isotropic on }\partial \Omega, 
\end{equation}
and, for $x \in \partial \Omega$, 
\begin{equation}\label{thm-reg-as3}
|\epsilon(x)-\hat{\epsilon}(x)| \geq \Lambda_1, \phantom{aaa}|\mu(x)-\hat{\mu}(x)|\geq \Lambda_1, \phantom{aaa}|\epsilon(x)/\mu(x)-\hat{\epsilon}(x)/\hat{\mu}(x)| \geq \Lambda_1.
\end{equation}
There exist two constants $k_0 \ge 1$ and $C > 0$  such that,  for $k \in \mC$ with $|\Im{(k^2)}| \ge \gamma |k|^2$ and $|k| \ge k_0$,  the conclusion of \Cref{thm-CN} holds for  $(J_e, J_m, \hJ_e, \hJ_m) \in [L^2(\Omega)]^{12}$. Moreover, for   $J_e,J_m,\hat{J}_e,\hat{J}_m \in [H^1(\Omega)]^3$ with $\dive J_e,\dive J_m,\dive \hat{J}_e,\dive \hat{J}_m \in H^1(\Omega)$ and $J_e \cdot \nu -\hat{J}_e \cdot \nu, J_m \cdot \nu -\hat{J}_m \cdot \nu \in H^{3/2}(\partial \Omega)$,  we have 
\begin{multline}\label{thm-reg-cl}
\| (E, H, \hE, \hH) \|_{H^2(\Omega)} + |k| \| (E, H, \hE, \hH) \|_{H^1(\Omega)} + 
 |k|^2 \,  \| (E, H, \hE, \hH) \|_{L^2(\Omega)}   \\[6pt]
 \le C |k| \| (J_e, J_m, \hJ_e, \hJ_m)\|_{L^2(\Omega)}
 + C  \| (J_e, J_m, \hJ_e, \hJ_m)\|_{H^1(\Omega)}  \\[6pt]
  + C  \| (\dive J_e,\dive  J_m, \dive \hJ_e,\dive  \hJ_m)\|_{L^2(\Omega)} +  \frac{C}{|k|}  \| (\dive J_e,\dive  J_m, \dive \hJ_e,\dive  \hJ_m)\|_{H^1(\Omega)}  \\[6pt]
 + C  \| (J_{e} \cdot \nu - \hJ_e \cdot \nu, J_m \cdot \nu -  \hJ_m \cdot \nu)\|_{H^{1/2}(\partial \Omega)} +  \frac{C}{|k|}  \| (J_{e} \cdot \nu - \hJ_e \cdot \nu, J_m \cdot \nu -  \hJ_m \cdot \nu)\|_{H^{3/2}(\partial \Omega)}, 
\end{multline}
for some positive constant $C$ depending only on $\Omega$, $\Lambda$, $\Lambda_1$, and $\gamma$. 
\end{theorem}

\begin{remark} \label{rem-iso} \rm The convention used in \eqref{H},  and in \eqref{thm-reg-as3}  are as follows.  A $3 \times 3$ matrix-valued function $M$ defined in a  subset $O \subset \mR^3$ is called isotropic at $x \in O$ if it is proportional  to the identity matrix at $x$, i.e., $M(x) = m I$ for some scalar $m = m(x)$ where $I$ denotes the $3 \times 3$ identity matrix. In this case, for the notational ease, we also denote $m(x)$ by $M(x)$. If $M$ is isotropic for $x \in O$, then $M$ is said to be isotropic in $O$. Condition \eqref{H} and \eqref{thm-reg-as3} are  understood under the convention $m(x) = M(x)$.
\end{remark}

Denote 
$$
\mR^3_+ = \Big\{x = (x_1, x_2, x_3) \in \mR^3; \; x_3 > 0 \Big\}
$$
and 
$$
\mR^3_0 = \Big\{x = (x_1, x_2, x_3) \in \mR^3; \; x_3 = 0 \Big\}. 
$$
One of the main ingredients of the proof of \Cref{thm-reg} is the following lemma, which is a variant of \cite[Corollary 3.1]{Cakoni-Ng21} (see also \Cref{rem-lem-HS}).

\begin{lemma}\label{lem-HS}  
Let $\gamma > 0$,  $k  \in \mC$ with $|\Im(k^2)| \ge \gamma |k|^2$,  and $|k| \ge 1$, and   let $\eps, \,  \mu,  \, \heps, \,  \hmu \in [C^1(\bar \mR^3_+)]^{3 \times 3}$ be symmetric, uniformly elliptic. Let $\Lambda \ge 1$ be such that
$$
\Lambda^{-1} \le \eps, \, \mu, \,  \heps, \,  \hmu  \le \Lambda \mbox{ in } B_1 \cap \mR^3_+  \quad  \mbox{ and } \quad \|(\eps, \mu, \heps, \hmu)  \|_{C^1(\mR^3_+ \cap B_1)} \le \Lambda. 
$$
Assume that $\eps(0), \, \heps(0), \, \mu(0), \hmu(0)$ are isotropic,  and for some $\Lambda_1 \ge 0$
$$
|\eps (0) -  \heps(0)| \ge \Lambda_1, \quad |\mu(0) - \hmu(0)| \ge \Lambda_1, \quad \mbox{ and } \quad |\eps(0)/ \mu (0) - \heps(0)/ \hmu(0)| \ge \Lambda_1. 
$$
Let $J_e, J_m, \hJ_e, \hJ_m \in [L^2(\mR^3_+)]^3$ and assume that  $(E, H, \hE, \hH) \in [L^2(\mR^3)]^{12}$ be a solution of the system\footnote{Here and in what follows $e_3=(0,0,1)$.}
\begin{equation}\label{lem-HS-sys}
\left\{\begin{array}{c}
\nabla \times E = k \mu H + J_e \mbox{ in } \mR^3_+, \\[6pt]
\nabla \times H = - k \eps E + J_m \mbox{ in } \mR^3_+, 
\end{array}\right. \quad \left\{\begin{array}{c}
\nabla \times \hE = k \hmu \hH + \hJ_e \mbox{ in } \mR^3_+, \\[6pt]
\nabla \times \hH = -  k \heps \hE + \hJ_m \mbox{ in } \mR^3_+, 
\end{array}\right.
\end{equation}
\begin{equation}\label{lem-HS-bdry}
(\hE - E) \times e_3 = 0 \mbox{ on } \mR^3_0, \quad \mbox{ and } \quad (\hH - H) \times e_3 = 0 \mbox{ on } \mR^3_0.  
\end{equation}
There exist $0 < r_0 < 1$ and $k_0 > 1$  depending only on $\gamma$, $\Lambda$, and $\Lambda_1$ such that if   the supports of $E, \, H, \, \hE, \, \hH$ are in $B_{r_0} \cap \overline{\mR^3_+}$, 
then, for $|k| > k_0$,  
\begin{enumerate}
\item[$i)$] 
\begin{equation}\label{lem-HS-CN-est1}
 |k| \,  \| (E, H, \hE, \hH) \|_{L^2(\mR^3_+)} 
 \le C   \| (J_e, J_m, \hJ_e, \hJ_m)\|_{L^2(\mR^3_+)}.  
\end{equation}

\item[$ii)$] if  $J_e, J_m, \hJ_e, \hJ_m \in H(\dive, \mR^3_+)$ and $J_{e, 3} - \hJ_{e, 3}, J_{m, 3} -  \hJ_{m, 3} \in H^{1/2}(\mR^3_0)$, then 
\begin{multline}\label{lem-HS-CN-est2}
\| (E, H, \hE, \hH) \|_{H^1(\mR^3_+)} + |k| \,  \| (E, H, \hE, \hH) \|_{L^2(\mR^3_+)} 
 \le C \Big(   \| (J_e, J_m, \hJ_e, \hJ_m)\|_{L^2(\mR^3_+)} \\[6pt] + \frac{1}{|k|} \| (\dive J_e, \dive J_m, \dive \hJ_e, \dive \hJ_m)\|_{L^2(\mR^3_+)} +  \frac{1}{|k|} \| (J_{e, 3} - \hJ_{e, 3}, J_{m, 3} -  \hJ_{m, 3})\|_{H^{1/2}(\mR^3_0)}  \Big).  
\end{multline}

\item[$iii)$] assume in addition that $\eps, \,  \mu,  \, \heps, \,  \hmu \in [C^2(\bar \mR^3_+)]^{3 \times 3}$ and 
\[\|(\eps, \mu, \heps, \hmu)  \|_{C^2( \mR^3_+ \cap B_1)} \le \Lambda.\]
Then, if $J_e, J_m, \hJ_e, \hJ_m \in [H^1(\mR^3_+)]^3$, $\dive J_e, \dive J_m, \dive \hJ_e, \dive \hJ_m \in H^1(\mR^3_+)$, and $J_{e, 3} - \hJ_{e, 3},  \, J_{m, 3} - \hJ_{m, 3} \in H^{3/2} (\mR^3_0)$, we have
\begin{align}\label{lem-HS-est}
 \| & (E, H, \hE, \hH)   \|_{H^2(\mR^3_+)} +  |k| \| (E, H, \hE, \hH)   \|_{H^1(\mR^3_+)}  + |k|^2 \,  \| (E, H, \hE, \hH) \|_{L^2(\mR^3_+)}  \\[6pt]
 \quad \le  & C |k| \| (J_e, J_m, \hJ_e, \hJ_m)\|_{L^2(\mR^3_+)} + C
 \| (J_e, J_m, \hJ_e, \hJ_m)\|_{H^1(\mR^3_+)} \nonumber \\[6pt]
 & +  C \| (\dive J_e, \dive J_m, \dive \hJ_e, \dive \hJ_m)\|_{L^2(\mR^3_+)} +   \frac{C}{|k|}  \| (\dive J_e, \dive J_m, \dive \hJ_e, \dive \hJ_m)\|_{H^1(\mR^3_+)} \nonumber \\[6pt]
& +  C \| (J_{e, 3} - \hJ_{e, 3}, J_{m, 3} - \hJ_{m, 3}) \|_{H^{1/2}(\mR^3_0)} + \frac{C}{|k|}  \| (J_{e, 3} - \hJ_{e, 3}, J_{m, 3} - \hJ_{m, 3}) \|_{H^{3/2}(\mR^3_0)}.   \nonumber
\end{align}
\end{enumerate}
Here $C$ denotes a positive constant depending only on $\gamma$,  $\Lambda$ and $\Lambda_1$. 
\end{lemma}

\begin{remark} \label{rem-lem-HS} \rm Parts $i)$ and $ii)$ are from \cite[Corollary 3.1]{Cakoni-Ng21},  which are restated here for the convenience of the reader. The new material is part $iii)$. 
\end{remark}

\begin{proof} We only prove $iii)$ (see \Cref{rem-lem-HS}). The idea of the proof  is as follows.  To derive \eqref{lem-HS-est}, we first differentiate the system with respect to $x_j$ for $j=1, 2$ and then derive the corresponding estimates for $(\partial_{x_j} E, \partial_{x_j} H, \partial_{x_j} \hE, \partial_{x_j} \hH)$ using $i)$ and $ii)$. After that,  we use the system of $(E, H)$ and $(\hE, \hH)$ to obtain similar estimates for $(\partial_{x_3} E, \partial_{x_3} H, \partial_{x_3} \hE, \partial_{x_3} \hH)$. This strategy is quite standard at least in the regularity theory of second elliptic equations, see e.g. \cite{Brezis-FA}. The main goal of the process is to keep track the dependence on $|k|$. The details are now given. 

Fix $k_0$ and $r_0$ such that $i)$ and $ii)$ hold.  By $ii)$, we have 
\begin{multline}\label{lem-HS-est1}
\| (E,  H, \hE, \hH) \|_{H^1(\mR^3_+)} + |k| 
\| (E,  H, \hE, \hH) \|_{L^2(\mR^3_+)}  \\[6pt]
 \le C  \left(  \|  (J_e, J_m, \hJ_e, \hJ_m) \|_{L^2(\mR^3_+)} + \frac{1}{|k|}  \| (\dive J_e, \dive J_m, \dive \hJ_e, \dive \hJ_m) \|_{L^2(\mR^3_+)} \right. \\[6pt]+ \left. \frac{1}{|k|}  \| (J_{e, 3} - \hJ_{e, 3}, J_{m, 3} - \hJ_{m, 3}) \|_{H^{1/2}(\mR^3_0)} \right).  
\end{multline}

Let $j =1, 2$.  Differentiating \eqref{lem-HS-sys} and \eqref{lem-HS-bdry} with respect to $x_j$, we obtain 
\begin{equation*} 
\left\{\begin{array}{c}
\nabla \times \partial_{x_j} E = k \mu \partial_{x_j}  H +   \bJ_e \mbox{ in } \mR^3_+, \\[6pt]
\nabla \times \partial_{x_j}  H = - k \eps \partial_{x_j}  E +   \bJ_m \mbox{ in } \mR^3_+, 
\end{array}\right. \quad \left\{\begin{array}{c}
\nabla \times \partial_{x_j}  \hE =  k \hmu \partial_{x_j}  \hH +   \hbJ_e \mbox{ in } \mR^3_+, \\[6pt]
\nabla \times \partial_{x_j}\hH = -  k \heps \partial_{x_j}   \hE + \hbJ_m \mbox{ in } \mR^3_+, 
\end{array}\right.
\end{equation*}
\begin{equation*}
(\partial_{x_j}  \hE - \partial_{x_j}  E) \times e_3 = 0 \mbox{ on } \mR^3_0, \quad \mbox{ and } \quad (\partial_{x_j}  \hH - \partial_{x_j}  H) \times e_3 = 0 \mbox{ on } \mR^3_0,   
\end{equation*}
where 
$$
\bJ_e = \partial_{x_j} J_e + k (\partial_{x_j} \mu) H, \quad 
\bJ_m = \partial_{x_j} J_m - k (\partial_{x_j} \eps) E, 
$$
$$ 
\hbJ_e = \partial_{x_j} \hJ_e + k (\partial_{x_j} \hmu) \hH, \quad 
\hbJ_m = \partial_{x_j} \hJ_m - k (\partial_{x_j} \heps) \hE. 
$$

Applying $ii)$ to  $(\partial_{x_j} E, \partial_{x_j} H, \partial_{x_j} \hE, \partial_{x_j} \hH)$, we deduce that 
\begin{align}\label{lem-HS-est2}
\| (\partial_{x_j} E,  \partial_{x_j} H, \partial_{x_j} \hE, \partial_{x_j} \hH) &  \|_{H^1(\mR^3_+)}   + |k| 
\| (\partial_{x_j} E,  \partial_{x_j} H, \partial_{x_j} \hE, \partial_{x_j} \hH) \|_{L^2(\mR^3_+)}  
 \le C (R_1 + R_2), 
\end{align}
where 
\begin{align}\label{lem-HS-R1}
R_1 = &   \|  (\partial_{x_j} J_e, \partial_{x_j} J_m, \partial_{x_j} \hJ_e, \partial_{x_j} \hJ_m) \|_{L^2(\mR^3_+)} \\[6pt]
 & + \frac{1}{|k|}  \| (\dive \partial_{x_j} J_e, \dive \partial_{x_j} J_m, \dive \partial_{x_j} \hJ_e, \dive \partial_{x_j} \hJ_m) \|_{L^2(\mR^3_+)} \nonumber \\[6pt]
 &+ \frac{1}{|k|}  \| (\partial_{x_j} J_{e, 3} - \partial_{x_j} \hJ_{e, 3}, \partial_{x_j} J_{m, 3} - \partial_{x_j} \hJ_{m, 3}) \|_{H^{1/2}(\mR^3_0)}, \nonumber  
\end{align}
and 
\begin{equation}\label{lem-HS-R2}
R_2 =    |k| \|  (E, H, \hE, \hH) \|_{L^2(\mR^3_+)} + \|  (E, H, \hE, \hH) \|_{H^1(\mR^3_+)} + \|  (E, H, \hE, \hH) \|_{H^{1/2}(\mR^3_0)}. 
\end{equation}
Combing \eqref{lem-HS-est1}, \eqref{lem-HS-R1},  and \eqref{lem-HS-R2}, we derive  from \eqref{lem-HS-est2}   that 
\begin{multline}\label{lem-HS-est3}
\| (\partial_{x_j} E,  \partial_{x_j} H, \partial_{x_j} \hE, \partial_{x_j} \hH)   \|_{H^1(\mR^3_+)}   + |k| 
\| (\partial_{x_j} E,  \partial_{x_j} H, \partial_{x_j} \hE, \partial_{x_j} \hH) \|_{L^2(\mR^3_+)} \\[6pt] 
 \le \mbox{the RHS of }  \eqref{lem-HS-est}. 
\end{multline}

On the other hand, from the system of $(E, H)$, we have, in $\mR^3_+$,  
\begin{multline}\label{lem-HS-relation}
\partial_{x_3}E_2=\partial_{x_2}E_3 - k(\mu H)_1-J_{e,1},\quad \partial_{x_3} E_1 = \partial_{x_1} E_3 + k (\mu  H)_2 +  J_{e, 2}\quad  \mbox{ and } \\[6pt] 
\partial_{x_3} \left ( \sum_{j=1}^3 \epsilon_{3j}E_j \right ) = -\sum_{\ell =1}^2 \sum_{j=1}^3 \partial_{x_\ell} \epsilon_{\ell j}E_j +\frac{1}{k}\dive (J_m) . 
\end{multline}
Combining \eqref{lem-HS-est1}, \eqref{lem-HS-est3},  and \eqref{lem-HS-relation}, and using the fact that $\eps_{33} \ge \Lambda^{-1}$, one has
\begin{equation}\label{lem-HS-est4}
\| E  \|_{H^2(\mR^3_+)} +  |k| \| E  \|_{H^1(\mR^3_+)}  + |k|^2 \,  \| E \|_{L^2(\mR^3_+)}  \le \mbox{the RHS of }  \eqref{lem-HS-est}. 
\end{equation}
Similarly, one obtains
\begin{equation}\label{lem-HS-est5}
\| (H, \hE, \hH)  \|_{H^2(\mR^3_+)} +  |k| \|  (H, \hE, \hH)  \|_{H^1(\mR^3_+)}  + |k|^2 \,  \|  (H, \hE, \hH)  \|_{L^2(\mR^3_+)}  \le \mbox{the RHS of }  \eqref{lem-HS-est}.
\end{equation}
The conclusion of \Cref{lem-HS} follows from \eqref{lem-HS-est3}, \eqref{lem-HS-est4}, and \eqref{lem-HS-est5}. 
\end{proof}

We are ready to give 
\begin{proof}[Proof of \Cref{thm-reg}] Let $K$ be a compact subset of $\Omega$. Fix $\varphi \in C^2_c(\Omega)$ such that $\varphi = 1$ in $K$.  Set 
$$
(E_\varphi, H_\varphi, \hE_\varphi, \hH_\varphi) = \varphi (E, H, \hE, \hH) \mbox{ in } \Omega.  
$$
From the system of $(E, H, \hE, \hH)$, we have 
\begin{equation}\label{sys-ITE-1}
\left\{\begin{array}{c}
\nabla \times E_\varphi = k \mu H_\varphi  + J_{\varphi,  e} \mbox{ in } \Omega, \\[6pt]
\nabla \times H_\varphi  = - k \eps E_\varphi  + J_{\varphi,m} \mbox{ in } \Omega, 
\end{array}\right. \quad \left\{\begin{array}{c}
\nabla \times \hE_\varphi  =  k \hmu \hH_\varphi  + \hJ_{\varphi, e} \mbox{ in } \Omega, \\[6pt]
\nabla \times \hH_\varphi  = -  k \heps \hE_\varphi  + \hJ_{\varphi, m} \mbox{ in } \Omega,  
\end{array}\right.
\end{equation}
\begin{equation}\label{bdry-ITE-1}
(\hE_\varphi - E_\varphi) \times \nu = 0 \mbox{ on } \partial \Omega, \quad \mbox{ and } \quad (\hH_\varphi - H_\varphi) \times \nu = 0 \mbox{ on } \partial \Omega. 
\end{equation}
Here, in $\Omega$,  
$$
J_{\varphi, e} = \nabla \varphi \times E + \varphi J_{e}, \; \; J_{\varphi, m} = \nabla \varphi \times H + \varphi J_{m}, \; \; \hJ_{\varphi, e} = \nabla \varphi \times \hE + \varphi \hJ_{e}, \; \;  \hJ_{\varphi, m} = \nabla \varphi \times \hH + \varphi \hJ_{m}. 
$$
Differentiating the system of $(E_\varphi, H_\varphi, \hE_\varphi, \hH_\varphi)$ with respect to $x_j$ ($1 \le j \le 3$) and applying \Cref{thm-CN}, we obtain, as in the proof of \Cref{lem-HS},  
\begin{multline*}
\| (E_\varphi, H_\varphi, \hE_\varphi, \hH_\varphi) \|_{H^2(\Omega)} 
 \le C |k| \| (J_{\varphi, e}, J_{\varphi, m}, \hJ_{\varphi, e}, \hJ_{\varphi, m})\|_{L^2(\Omega)}
 + C  \| (J_{\varphi, e}, J_{\varphi, m}, \hJ_{\varphi, e}, \hJ_{\varphi, m})\|_{H^1(\Omega)}  \\[6pt]
  + C  \| (\dive J_{\varphi, e}, \dive J_{\varphi, m}, \dive \hJ_{\varphi, e}, \dive \hJ_{\varphi, m}) \|_{L^2(\Omega)} +  \frac{C}{|k|}  \|(\dive J_{\varphi, e}, \dive J_{\varphi, m}, \dive \hJ_{\varphi, e}, \dive \hJ_{\varphi, m})\|_{H^1(\Omega)}. 
\end{multline*} 
This implies
\begin{multline}\label{thm-reg-p11}
\| (E_\varphi, H_\varphi, \hE_\varphi, \hH_\varphi) \|_{H^2(\Omega)} 
 \le C |k| \| (J_e, J_m, \hJ_e, \hJ_m)\|_{L^2(\Omega)}
 + C  \| (J_e, J_m, \hJ_e, \hJ_m)\|_{H^1(\Omega)}  \\[6pt]
  + C  \| (\dive J_e, \dive J_m, \dive \hJ_e, \dive \hJ_m)\|_{L^2(\Omega)} +  \frac{C}{|k|}  \|(\dive J_e, \dive J_m, \dive  \hJ_e, \dive \hJ_m)\|_{H^1(\Omega)}\\[6pt]
  + C |k| \| (E, H, \hE, \hH)\|_{L^2(\Omega)}
 + C  \| (E, H, \hE, \hH) \|_{H^1(\Omega)}. 
\end{multline} 
Applying \Cref{thm-CN} again, we derive from \eqref{thm-reg-p11} that 
\begin{multline}\label{thm-reg-p1}
\| (E_\varphi, H_\varphi, \hE_\varphi, \hH_\varphi) \|_{H^2(\Omega)} + 
|k| \| (E_\varphi, H_\varphi, \hE_\varphi, \hH_\varphi) \|_{H^1(\Omega)}
\\[6pt] 
+ |k|^2 \| (E_\varphi, H_\varphi, \hE_\varphi, \hH_\varphi) \|_{L^2(\Omega)} 
 \le \mbox{the RHS of } \eqref{thm-reg-cl}. 
\end{multline} 

The conclusion of \Cref{thm-reg} now follows from \eqref{thm-reg-p1} and \Cref{lem-HS} via local charts. The proof is complete. 
\end{proof}

\section{Completeness of the generalized eigenfunctions - Proof of \Cref{thm1}} \label{sect-thm1}

To establish the completeness of the generalized eigenfunctions, we use \Cref{thm-reg} and apply the theory of Hilbert-Schmidt operators. To this end, we first recall 

 \begin{definition} Let $H$ be a separable Hilbert space and let $(\phi_k)_{k=1}^\infty$ be an orthogonal basis. A bounded linear operator $\mathbf{T}: H \to H$ is Hilbert- Schmidt  if its finite double norm 
 \[
 \vvvert
 \mathbf{T}  \vvvert: =\left(\sum_{k=1}^{\infty} \| \mathbf{T}(\phi_k)\|_{H}^2\right)^{1/2} < + \infty. 
 \]
 \end{definition}

\begin{remark} \rm The definition of $ \vvvert {\mathbf T}  \vvvert$ does not depend on the choice of $(\phi_k)$, see e.g. \cite[Chapter 12]{Agmon}. 
\end{remark}

Using \Cref{thm-reg}, we can establish the following result.

\begin{proposition} \label{pro-HS} Assume that  $\eps, \, \mu, \, \heps, \, \hat{\mu}  \in [C^2(\bar \Omega)]^{3 \times 3}$ and \eqref{H} holds,  and  let $\gamma > 0$. Let $k_0 \ge 1$ and $C > 0$  be constants such that   for $k \in \mC$ with $|\Im{(k^2)}| \ge \gamma |k|^2$ and $|k| \ge k_0$, the conclusions  of \Cref{thm-reg} hold. Then, for such a complex number $k$,  
\begin{multline} \label{pro-HS-state}
\|\cT_k^2 (\bcJ) \|_{H^2(\Omega)} + |k| \|\cT_k^2 (\bcJ) \|_{H^1(\Omega)} + |k|^2 \|\cT_k^2 (\bcJ) \|_{L^2(\Omega)} \\[6pt]
\le C \| \cJ  \|_{L^2(\Omega)} \quad \forall \cJ= (\cJ_e, \cJ_m, \hcJ_e, \hcJ_m) \in \bH(\Omega).  
\end{multline}
Consequently, 
\begin{itemize}
\item[i)] $\cT_k^2$ is a Hilbert-Schmidt operator defined on $\bH(\Omega)$; moreover, 
\begin{equation} \label{pro-HS-T2}
\vvvert \cT_k^2 \vvvert \le \frac{C}{|k|^{1/2}}, 
\end{equation}
for some positive constant $C$, independent of $k$.

\item[ii)]  For $\theta \in \mR$ with $|\Im (e^{2 i \theta})| > 0$, $e^{i \theta}$  is a direction of minimal growth of the modified resolvent of $\cT_k^2$.  
\end{itemize}
\end{proposition}

For the convenience of the reader, we recall briefly here some notions associated to the concept of the minimal growth. Let $A$ be a continuous, linear transformation from a Hilbert space $H$ into itself.  
The modified resolvent set $\rho_m(A)$ of $A$ is the set of all $\lambda \in \mC \setminus \{0 \}$ such that $I - \lambda A$ is bijective (and continuous).  If $\lambda \in \rho_m(A)$, then the map $A_\lambda: = A (I - \lambda A)^{-1}$ is the modified resolvent of $A$ (see \cite[Definition 12.3]{Agmon}). For $\theta \in \mR$,   $e^{i \theta}$ is a direction of minimal growth of the modified resolvent of $A$ if for some $a>0$, the following two facts hold for all $r>a$:  $i)$ $r e^{i \theta}$ is in the modified resolvent set $\rho_m(A)$ of $A$ and $ii)$ $\| A_{r e^{i \theta}} \| \le C/ r$ (see \cite[Definition 12.6]{Agmon}).

\medskip 
Another key ingredient of  the proof of \Cref{thm1} is: 

\begin{proposition} \label{pro-density} Assume that  $\eps, \, \mu, \, \heps, \, \hat{\mu}  \in [C^2(\bar \Omega)]^{3 \times 3}$ and \eqref{H} holds.  Let $k \in \C$ be such that the conclusion of \Cref{thm-reg} holds. We have
\[
\overline{\cT_k(\bm{H}(\Omega))}^{L^2(\Omega)} = \bm{H}(\Omega).
\]
\end{proposition}

The rest of this section containing three subsections  is organized as follows.  In the first  subsection, we give the proof of \Cref{pro-HS}.  The proofs of \Cref{pro-density} and \Cref{thm1} are given in the last two subsections, respectively.

\subsection{Proof of \Cref{pro-HS}}

We first state and prove a simple useful lemma used in the proof of \Cref{pro-HS}. 

\begin{lemma}\label{lem-Tmod} Let $k, \, s  \in \mC$ be such that $\cT_k, \cT_{k+s} : \bH(\Omega) \to \bH(\Omega)$ are bounded.  We have 
\begin{enumerate}

\item[i)] If $\cT_k$ is compact, then $s \in \rho_m(\cT_k). $

\item[ii)] Assume  $s \in \rho_m(\cT_k)$. Then 
\begin{equation}
 \cT_k (I - s \cT_k)^{-1} = (I - s \cT_k)^{-1}  \cT_k  = \cT_{k + s}.
\end{equation}

\end{enumerate}
\end{lemma}

\begin{proof}[Proof of \Cref{lem-Tmod}] We begin with assertion $i)$. Since $\cT_k$ is compact, it suffices to prove that $I - s \cT_k$ is injective. Indeed, let $(E, H, \hE, \hH) \in \bH(\Omega)$ be a solution of the equation $I - s \cT_k = 0$. One can check that   $(E, H, \hE, \hH)  = \cT_{k+s} (0) = 0$. Assertion $i)$ follows.

We next establish $ii)$. Let $\cJ = (\cJ_e, \cJ_m, \hcJ_e, \hcJ_m) \in \bH(\Omega)$ be arbitrary. Set 
\begin{equation}\label{lem-Tmod-1}
(E, H, \hE, \hH) = \cT_{k+s} (\cJ), 
\end{equation}
\begin{equation}\label{lem-Tmod-2}
\cJ^1 = (\cJ_e^1, \cJ_m^1, \hcJ_e^1, \hcJ_m^1)= (I - s \cT_k)^{-1} (\cJ), 
\end{equation}
\begin{equation}\label{lem-Tmod-3}
(E^1, H^1, \hE^1, \hH^1) = \cT_{k} (\cJ^1). 
\end{equation}
We claim that 
$$
(E^1, H^1, \hE^1, \hH^1)   = (E, H, \hE, \hH), 
$$
which implies $\cT_k (I - s \cT_k)^{-1} = \cT_{k+s}$ since $\cJ$ is arbitrary. 

To prove the claim, we will show that $(E^1, H^1, \hE^1, \hH^1)$ and $(E, H, \hE, \hH)$ satisfy the same Cauchy problem.  We have
\begin{equation}\label{lem-Tmod-4}
\nabla \times E \mathop{=}^{\eqref{lem-Tmod-1}} (k  + s) \mu H + \mu \cJ_m, 
\end{equation}
$$
\nabla \times E^1 \mathop{=}^{\eqref{lem-Tmod-3}} k \mu H^1 + \mu \cJ_m^1, 
$$
$$
\cJ^1 - \cJ \mathop{=}^{\eqref{lem-Tmod-2}} s \cT_k (\cJ_1) \mathop{=}^{\eqref{lem-Tmod-3}} s (E^1, H^1, \hE^1, \hH^1). 
$$
This implies 
$$
\nabla \times E^1 = (k + s) \mu H^1 + \mu \cJ_m, 
$$
(compare with \eqref{lem-Tmod-4}). 
Similarly, we can derive that  $(E^1, H^1, \hE^1, \hH^1)$ and $(E, H, \hE, \hH)$ satisfy the same system  since it is clear that, on $\partial \Omega$,  
$$
(\hE_1 - E_1)\times \nu =  (\hH_1 - H_1)\times \nu = 0 = (\hE - E)\times \nu =  (\hH - H)\times \nu. 
$$
The claim is proved. 

Since 
$$
(I - s \cT_k) (I - s \cT_k)^{-1} = I =  (I - s \cT_k)^{-1} (I - s \cT_k), 
$$
and $s \neq 0$ by the definition of $\rho_m(\cT_k)$, we obtain 
$$
 \cT_k (I - s \cT_k)^{-1} =   (I - s \cT_k)^{-1}  \cT_k.  
$$
The proof is complete. \end{proof}

We are ready to give 

\begin{proof}[Proof of \Cref{pro-HS}] Assertion \eqref{pro-HS-state} is just a consequence of \Cref{thm-CN} and \Cref{thm-reg}.  As a consequence of \eqref{pro-HS} and Gagliardo-Nirenberg's inequality see \cite{Gagliardo59,Nirenberg59}, we derive, for $\cJ \in \bH(\Omega)$, that $\cT_k^2(\cJ) \in [C(\bar \Omega)]^{12}$, and 
$$
\| \cT_k^2(\cJ) \|_{L^\infty(\Omega)}
\le C \| \cT_k^2(\cJ) \|_{H^2(\Omega)}^\frac{3}{4} \| \cT_k^2(\cJ) \|_{L^2(\Omega)}^\frac{1}{4} \le \frac{C}{|k|^{1/2}} \| \cJ \|_{L^2(\Omega)}. 
$$
It follows from the theory of Hilbert-Schmidt operators, see e.g. \cite[Lemma 3]{MinhHung2} \footnote{In \cite[Lemma 3]{MinhHung2}, the statement is on $[L^2(\Omega)]^m$ for some $m \ge 1$, nevertheless, the proof also gives the result for $\bH(\Omega)$ since $\bH(\Omega)$ is equipped with the $[L^2(\Omega)]^{12}$-norm.}, that  $\cT_k^2$ is a Hilbert-Schmidt operator defined on $\bH(\Omega)$ and 
$$
 \vvvert  \cT_k^2 \vvvert \le   \frac{C}{|k|^{1/2}}. 
$$ 

We next check the assertion on the minimal growth of the modified resolvent of $\cT_k$. 
We have 
$$
\lim_{r \to + \infty} |\Im \big( (k  + r e^{i \theta})^2 \big)| / |k  + r e^{i \theta}|^2  =  |\Im (e^{2 i \theta}) | \ge 2 \gamma,  
$$
for some $\gamma > 0$.  It follows, for $a$ large enough,  that $k + r e^{i \theta}$ satisfies the conclusion of \Cref{thm-reg} for $r > a$. 
On the other hand, 
let $(E, H, \hE, \hH) \in \bH(\Omega)$.  We first note that, for $s \in \mC$,  
$$
(I - s \cT_k) (E, H, \hE, \hH) = 0 \mbox{ if and only if }  (E, H, \hE, \hH) = \cT_{k + s} ( 0) = 0, 
$$
provided that $\cT_{k+s}$ is well-defined.  Since $\cT_k$ is compact, it follows that $ r e^{i \theta} \in \rho_m(\cT_k)$ for $r > a$. By \Cref{lem-Tmod}, we also have, with $s = r e^{i \theta}$,  
$$
 \cT_k (I - s \cT_k)^{-1} = (I - s \cT_k)^{-1}  \cT_k  = \cT_{k + s}. 
$$
Let $s_1= i r^{1/2} e^{i \theta/2}$ and $s_2 = - i r^{1/2} e^{i \theta/2}$. Thus $(t- s_1)(t-s_2) = t^2 - s$ for $t \in \mC$. One then can check that  
\begin{equation*}
\cT_k^2 (I - s \cT_k^2)^{-1} = \cT_k^2 (I - s_1 \cT_k)^{-1}(I - s_2 \cT_k)^{-1}  \mathop{=}^{\Cref{lem-Tmod}}  \cT_k(I - s_1 \cT_k)^{-1} \cT_k (I - s_2 \cT_k)^{-1}  = \cT_{k+ s_1} \cT_{k+s_2}. 
\end{equation*}
It follows from \Cref{thm-CN} that 
\begin{align*}
\| \cT_k^2 (I - s \cT_k^2)^{-1}  \|_{\bH(\Omega) \to \bH(\Omega)} &=  \| \cT_{k+ s_1} \cT_{k+s_2} \|_{\bH(\Omega) \to \bH(\Omega)}\\[6pt] 
&\le \| \cT_{k+ s_1} \|_{\bH(\Omega) \to \bH(\Omega)} \|\cT_{k+s_2} \|_{\bH(\Omega) \to \bH(\Omega)} \\[6pt]
&\le C \frac{1}{|s_1|}\frac{1}{|s_2|} =  \frac{C}{|s|}.  \nonumber
\end{align*}
The assertion on the minimal growth of the modified resolvent  of $\cT_k^2$ follows. 
\end{proof}

\subsection{Proof of  \Cref{pro-density}} 

We first state and prove the following technical result which is used in the proof of \Cref{pro-density}. 

\begin{lemma} \label{lem-density1} 
Let  $\M \in [C^1(\bar{\Omega})]^{3 \times 3}$ be symmetric and uniformly elliptic. Let $U \in [H^1(\Omega)]^3$ be such that $\dive (\M U) = 0$ in $\Omega$. There exists a  sequence $(U_n)_n \subset [H^1(\Omega)]^3$ such that 
\begin{equation}
\dive (\M U_n) = 0 \mbox{ in } \Omega,
\end{equation}
\begin{equation}
\M U_n \cdot \nu = \M U \cdot \nu \mbox{ on } \partial \Omega, \quad U_n \times \nu = 0 \mbox{ on } \partial \Omega, 
\end{equation}
 and
\begin{equation}
U_n \to U  \text{ in } [L^2(\Omega)]^3 \mbox{ as } n \to + \infty. 
\end{equation}
\end{lemma}

\begin{proof}[Proof of \Cref{lem-density1}] 
Since $\Omega$ is connected, $U \in [H^1(\Omega)]^3$ and $\dive (\M U) = 0$ in $\Omega$, by  \cite[lemma 2.2]{GR86}, there exists $\widetilde{V} \in [H^1(\Omega)]^3$ such that 
\begin{equation}\label{lem-density1-tVn}
\dive \widetilde{V}  = 0 \text{ in }\Omega \: \: \: \text{ and }\: \: \: \widetilde{V} = \frac{ \mathcal{M} U \cdot  \nu  }{ \M \nu \cdot  \nu } \M \nu  \text{ on }\partial \Omega.
\end{equation}
Set $V = \mathcal{M}^{-1}\widetilde{V}$ in $\Omega$. One can easily check from the definition of $V$ and \eqref{lem-density1-tVn} that
\begin{equation}\label{lem-density1-Vn}
\dive (\mathcal{M}V)  = 0 \text{ in }\Omega, \: \: \: \M V \cdot \nu  =   \M U \cdot  \nu   \text{ on }\partial \Omega \: \: \:  \text{ and }\: \: \: V \times \nu = 0  \text{ on }\partial \Omega.
\end{equation}
Set $\widetilde{U}= U-V$ in $\Omega$.  Since $\dive (\M U) = 0$ in $\Omega$, we derive from  \eqref{lem-density1-Vn} that $\dive (\M \widetilde{U})=0$ in $\Omega$ and $ \M \widetilde{U} \cdot \nu =0$ on $\partial \Omega$.  It follows from \cite[Theorem 2.8]{GR86} that there exists a sequence $(\widetilde{U}_n)_n \subset [C^1_c(\Omega)]^3$ such that 
\begin{equation}\label{lem-density1-U1n-1n}
\dive (\widetilde{U}_n)  = 0 \text{ in }\Omega 
\end{equation}
and 
\begin{equation}\label{lem-density1-U1n-2n}
\widetilde{U}_n \to \M \widetilde{U} \mbox{ in } [L^2(\Omega)]^3 \mbox{ as } n \to + \infty.
\end{equation}
Set
\[
U_n = \M^{-1} \widetilde{U}_n +  V.
\]
We claim that the sequence $(U_n)_n$ has the required properties. Indeed,
\[
\dive (\M U_n) =\dive (\widetilde{U}_n) + \dive (\M V) \mathop{=}^{\eqref{lem-density1-Vn},\eqref{lem-density1-U1n-1n}}  0 \text{ in }\Omega 
\]
and, since $\widetilde{U}_n \in [C_c^1(\Omega)]^3$, we also have
\[
\M U_n \cdot \nu =  \widetilde{U}_n \cdot \nu +  \M V \cdot \nu \mathop{=}^{\eqref{lem-density1-Vn}} \M U \cdot  \nu \text{ on }\partial \Omega, 
\]
and
\[
U_n \times \nu = \M^{-1} \widetilde{U}_n \times \nu  +  V \times \nu \mathop{=}^{\eqref{lem-density1-Vn}}  0\text{ on }\partial \Omega .
\]
Moreover, since $V \in [H^1(\Omega)]^3$, it follows that $U_n \in [H^1(\Omega)]^3$ and by \eqref{lem-density1-U1n-2n} we obtain 
\[
U_n \to \widetilde{U} + V = U \mbox{ in } [L^2(\Omega)]^3 \mbox{ as } n \to + \infty.
\]
The proof is complete.
\end{proof}

We are ready to give

\begin{proof}[Proof of  \Cref{pro-density}]
Since $\cT_k$ is a map from $\bH(\Omega)$ into $\bH(\Omega)$, it suffices to prove the following two facts
\begin{equation}\label{pro-density-part1}
[H^1(\Omega)]^{12} \cap \bH(\Omega) \subset \overline{\cT_k \big( \bH(\Omega) \big)}^{L^2(\Omega)},  
\end{equation}
and 
\begin{equation}\label{pro-density-part2}
[H^1(\Omega)]^{12} \cap \bH(\Omega) \mbox{ is dense in } \bH(\Omega) \mbox{ with respect to } [L^2(\Omega)]^{12}-\mbox{norm}.  
\end{equation}
These will be proved in Steps 1 and 2 below. 

{\it Step 1: Proof of \eqref{pro-density-part1}}.  Let $(E,H,\hat{E},\hat{H}) \in [H^1(\Omega)]^{12} \cap \bH(\Omega)$. By applying  \Cref{lem-density1} with $(\M, U)$ equal to $(\eps,  E)$, $(\mu,  H)$, $(\heps,  \hE)$, and $(\hmu,  \hH)$,  there exists a sequence $\big((E^n,H^n,\hat{E}^n,\hat{H}^n) \big)_n \subset [H^1(\Omega)]^{12} \cap \bH(\Omega)$ such that
\begin{equation}\label{pro-density-En1}
E^n \times \nu = H^n \times \nu = \hat{E}^n \times \nu = \hat{H}^n \times \nu = 0 \mbox{ on } \partial \Omega,
\end{equation}
and
\begin{equation}\label{pro-density-En2}
(E^n,H^n,\hat{E}^n,\hat{H}^n) \to (E,H,\hat{E},\hat{H})  \mbox{ in } [L^2(\Omega)]^{12} \mbox{ as } n \to + \infty. 
\end{equation}

Set, in $\Omega$,  
\begin{equation}\label{pro-density-J}
J_e^n = \nabla \times E^n - k \mu H^n, \quad  J_m^n =  \nabla \times H^n +  k \eps H^n, 
\end{equation}
\begin{equation}\label{pro-density-hJ}
\hJ_e^n = \nabla \times \hE^n - k \hmu \hH^n, \quad \hJ_m^n =  \nabla \times \hH^n +  k \heps \hH^n,  
\end{equation}
and  define $(\cJ_e^n, \cJ_m^n, \hcJ_e^n, \hcJ_m^n)$ in $\Omega$ via  $(J_e^n, J_m^n, \hJ_e^n, \hJ_m^n) = (\mu \cJ_m^n, - \eps \cJ_e^n, \hmu \hcJ_m^n, -  \heps \hcJ_e^n)$.  

It follows that \eqref{e11p} holds with $(E, H, \hE, \hH)$ and $(J_e, J_m, \hJ_e, \hJ_m)$ replaced by $(E^n, H^n, \hE^n, \hH^n)$ and $(J_e^n, J_m^n, \hJ_e^n, \hJ_m^n)$.  Since $(E^n,H^n,\hat{E}^n,\hat{H}^n) \in \bH(\Omega)$, it follows that 
\begin{equation}\label{pro-density-S1-p1}
\dive J_e^n  = \dive   J_m^n = \dive  \hJ_e^n = \dive \hJ_m^n  =0 \mbox{ in } \Omega. 
\end{equation}

On the other hand,  from \eqref{pro-density-J} and \eqref{pro-density-hJ}, we have, on $\partial \Omega$, 
\begin{equation*}
(\hJ_e^n - J_e^n ) \cdot \nu = (\nabla \times \hE^n - \nabla \times \hE^n) \cdot \nu - k  (\hmu \hH^n - \mu H^n) \cdot \nu.   
\end{equation*}
This implies 
\begin{equation}\label{pro-density-S1-p2}
(\hJ_e^n - J_e^n ) \cdot \nu  = 0 \mbox{ on } \partial \Omega, 
\end{equation}
since $(\hmu \hH^n - \mu H^n) \cdot \nu  = 0$ on $\partial \Omega$ and $\dive_{\partial \Omega} \Big( (\hE^n - E^n) \times \nu \Big) = 0$ on $\partial \Omega$ by \eqref{pro-density-En1}. Similarly, we have 
\begin{equation}\label{pro-density-S1-p3}
 (\hJ_m^n - J_m^n ) \cdot \nu = 0 \mbox{ on } \partial \Omega. 
\end{equation}

Combining \eqref{pro-density-S1-p1}, \eqref{pro-density-S1-p2}, and \eqref{pro-density-S1-p3} yields that $(\cJ_e^n, \cJ_m^n, \hcJ_e^n, \hcJ_m^n) \in \bH(\Omega)$. Consequently,
$$
(E^n, H^n, \hE^n, \hH^n) \in \cT_k \big( \bH(\Omega) \big). 
$$
The conclusion of Step 1 now follows from \eqref{pro-density-En2}.

\medskip 
{\it Step 2: Proof of \eqref{pro-density-part2}.} Fix 
$(E, \, H, \, \hat{E}, \, \hat{H})  \in \bH(\Omega)$ arbitrary.  There exist sequences $(\cE^n)_n, \, (\cH^n)_n \subset [H^2(\Omega)]^{3}$ such that
\begin{equation} \label{pro-density-part2-En1}
(\eps \cE^n, \, \mu \cH^n) \to (\eps E,  \, \mu H) \mbox{ in } [H(\dive, \Omega)]^{2}.  
\end{equation}
Since 
$$
\dive (\heps \hE - \eps E) = \dive (\hmu \hH - \mu H) = 0 \mbox{ in } \Omega \quad \mbox{ and } \quad 
(\heps \hE - \eps E) \cdot  \nu =  (\hmu \hH - \mu H) \cdot \nu = 0 \mbox{ on } \partial  \Omega, 
$$
by \cite[Theorem 2.8]{GR86}, there exist sequences  $(U^n_e)_n, \, (U_m^n)_n \subset [H^2(\Omega)]^3$ such that 
\begin{equation}\label{pro-density-U}
\dive U^n_e = \dive  U_m^n = 0 \mbox{ in } \Omega, 
\end{equation}
and
\begin{equation} \label{pro-density-part2-En2}
(U^n_e,  \, U_m^n) \to (\heps \hE - \eps E, \hmu \hH - \mu H) \mbox{ in } [L^2(\Omega)]^6  \mbox{ as } n \to + \infty. 
\end{equation} 

Define  $\hcE^n, \, \hcH^n \in [L^2(\Omega)]^3$ via 
\begin{equation}\label{pro-density-def-hEH}
\heps \hcE^n = U_e^n + \eps \cE^n \mbox{ in } \Omega \quad \mbox{ and } \quad \hmu \hcH^n = U_m^n + \hmu \cH^n \mbox{ in } \Omega. 
\end{equation}
From \eqref{pro-density-part2-En1}, \eqref{pro-density-U}, and \eqref{pro-density-part2-En2}, we have 
\begin{equation}\label{pro-density-part2-En3}
(\heps \hcE^n, \, \hmu \hcH^n) \to (\heps \hE,  \, \hmu \hH) \mbox{ in } [H(\dive, \Omega)]^{2}. 
\end{equation}

Using \eqref{pro-density-part2-En1} and \eqref{pro-density-part2-En3}, we derive from the trace theory that, as $n \to + \infty$,  
\begin{equation}\label{pro-density-trace}
(\eps \cE^n  - \eps E) \cdot \nu, \; 
(\mu \cH^n  - \mu H) \cdot \nu, \; 
(\heps \hcE^n  - \heps \hE) \cdot \nu, \; 
(\hmu \hcH^n  - \hmu \hH) \cdot \nu \to 0  \mbox{ in } H^{-1/2} (\partial \Omega). 
\end{equation}
Since $(\heps \hE  - \eps E) \cdot \nu =  
(\hmu \hH  - \mu H) \cdot \nu = 0 $ on $\partial \Omega$,  we obtain 
\begin{equation}\label{pro-density-trace1}
(\heps \hcE^n  - \eps \cE^n) \cdot \nu, \; 
(\hmu \hcH^n  - \hmu \cH^n) \cdot \nu \to 0, \; 
 \mbox{ in } H^{-1/2} (\partial \Omega) \mbox{ as } n \to + \infty. 
\end{equation}

Set 
\begin{equation}\label{pro-density-alpha}
\alpha_e^n = \frac{1}{|\partial \Omega|} \int_{\partial \Omega} \eps \cE^n \cdot \nu \quad \mbox{ and } \quad
\alpha_m^n = \frac{1}{|\partial \Omega|} \int_{\partial \Omega} \mu \cH^n \cdot \nu, 
\end{equation}
where $|\partial \Omega|$ denotes the 2-Hausdorff measure of $\partial \Omega$.  We derive  that 
\begin{equation}\label{pro-density-alpha-pro1}
\lim_{n \to + \infty} \alpha_e^n \mathop{=}^{\eqref{pro-density-trace}}\frac{1}{|\partial \Omega|} \int_{\partial \Omega} \eps E \cdot \nu = \frac{1}{|\partial \Omega|} \int_{ \Omega} \dive (\eps E ) = 0. 
\end{equation}
Similarly, we obtain 
\begin{equation}\label{pro-density-alpha-pro2}
\lim_{n \to + \infty} \alpha_m^n = 0. 
\end{equation}

Denote 
$$
H^1_{\sharp} (\Omega) = \Big\{u \in H^1(\Omega) : \int_\Omega u = 0 \Big\}.
$$
Let $\xi_e^n, \, \xi_m^n, \,  \hxi_e^n, \, \hxi_m^n \in H^1_{\sharp} (\Omega)$ be a  solution of 
\begin{equation}\label{pro-density-def-xi}
\left \{
\begin{array}{cl}
-\dive (\eps \nabla  \xi_e^n)  =  - \dive (\eps \cE^n) &  \text{ in }\Omega , \\[6pt]
\eps \nabla  \xi_e^n \cdot \nu  = \alpha_e^n & \text{ on }\partial \Omega ,
\end{array}
\right .
\quad
\left \{
\begin{array}{cl}
-\dive (\mu \nabla  \xi_m^n)  =  - \dive (\mu \cH^n)  &  \text{ in }\Omega , \\[6pt]
\mu \nabla  \xi_m^n \cdot \nu  = \alpha_m^n  & \text{ on }\partial \Omega ,
\end{array}
\right .
\end{equation}
\begin{equation}\label{pro-density-def-jxi}
\left \{
\begin{array}{cl}
-\dive (\heps \nabla  \hxi_e^n)  =  - \dive (\heps \hcE^n) &  \text{ in }\Omega , \\[6pt]
\heps \nabla  \hxi_e^n \cdot \nu  =  (\heps \hcE^n -  \eps \cE^n) \cdot \nu + \alpha_e^n  & \text{ on }\partial \Omega ,
\end{array}
\right .
\quad
\left \{
\begin{array}{cl}
-\dive (\hmu\nabla  \hxi_m^n)  =  - \dive (\hmu \hcH^n)  &  \text{ in }\Omega , \\[6pt]
\hmu \nabla  \hxi_m^n \cdot \nu  =  (\hmu \hcH^n  - \mu \cH^n) \cdot \nu + \alpha_m^n & \text{ on }\partial \Omega. 
\end{array}
\right .
\end{equation}

By the definition of $\alpha_e^n$ and $\alpha_m^n$ \eqref{pro-density-alpha}, we have 
\begin{equation}\label{pro-density-coucou}
\int_{\Omega} \dive (\eps \cE^n) = \int_{\partial \Omega} \alpha_e^n \quad \mbox{ and } \quad \int_{\Omega}  \dive (\mu \cH^n) = \int_{\partial \Omega} \alpha_m^n.  
\end{equation}
It follows that $\xi_e^n$ and $\xi_m^n$ are well-defined and uniquely determined.  We also have 
\[
\int_{\Omega} \dive (\heps \hcE^n)  - \int_{\partial \Omega}  \Big(  (\heps \hcE^n -  \eps \cE^n) \cdot \nu + \alpha_e^n \Big) \mathop{=}^{\eqref{pro-density-alpha}}0
\]
and 
$$
\int_{\Omega} \dive (\hmu \hcH^n)  - \int_{\partial \Omega}  \Big( (\hmu \hcH^n  - \mu \cH^n) \cdot \nu + \alpha_m^n  \Big) \mathop{=}^{\eqref{pro-density-alpha}} 0. 
$$
Hence $\hxi_e^n$ and $\hxi_m^n$ are well-defined and uniquely determined as well. From the regularity theory of elliptic equations it follows that
\begin{equation}
\label{pro-density-reg}
(\xi_e^n,\xi_m^n,\hat{\xi}^n_e,\hat{\xi}^n_m) \in [H^2(\Omega)]^4.
\end{equation}

Using \eqref{pro-density-trace1}, \eqref{pro-density-alpha-pro1}, and  \eqref{pro-density-alpha-pro2}, we derive that 
\begin{equation}\label{pro-density-xi1}
\xi_e^n, \, \xi_m^n, \, \hxi_e^n, \, \hxi_m^n \to  0 \mbox{ in } H^{1} (\Omega) \mbox{ as } n \to + \infty. 
\end{equation}

Set
\begin{equation}
(E_n, H_n, \hat{E}_n, \hat{H}_n) = (\cE^n   - \nabla \xi_e^n, \,  \cH^n -  \nabla \xi_m^n, \,  \hcE^n - \nabla \hxi_e^n, \, \hcH^n - \nabla \hxi_m^n) \mbox{ in } \Omega. 
\end{equation}
We have, by \eqref{pro-density-part2-En1}, and \eqref{pro-density-part2-En3}, and \eqref{pro-density-xi1}, 
\begin{equation}\label{pro-density-S2-cl1}
(E_n, H_n, \hat{E}_n, \hat{H}_n)  \to (E, H, \hE, \hH)\mbox{ in } [L^2(\Omega)]^{12}. 
\end{equation}
From the definition of $\xi_e^n, \, \xi_m^n, \,  \hxi_n^n, \,  \hxi_m^n$, we have  
\begin{equation}\label{pro-density-S2-cl2} 
\dive (\eps E^n)  = \dive( \heps \hE^n) = \dive (\mu H^n) = \dive (\hmu \hH^n)  = 0 \mbox{ in } \Omega,  
\end{equation}
and, on $\partial \Omega$,  
\begin{equation}\label{pro-density-S2-cl3} 
(\heps \hE^n  - \eps E^n) \cdot \nu = 
(\hmu \hH^n  - \mu H^n) \cdot \nu  =  0  \mbox{ on } \partial \Omega. 
\end{equation}

Combining \eqref{pro-density-reg}, \eqref{pro-density-S2-cl2}, and \eqref{pro-density-S2-cl3} yields 
\begin{equation}\label{pro-density-part2-p4}
(E^n, H^n, \hE^n, \hH^n) \in \bH(\Omega) \cap [H^1(\Omega)]^{12}.
\end{equation} 
The conclusion of Step 2 now follows from \eqref{pro-density-S2-cl1}. 

\medskip 
The proof is complete. 
\end{proof}

\begin{remark} \label{rem-Form} \rm  One can rewrite   \eqref{e11} and \eqref{e12} under the following form 
\begin{equation}\label{rem-Form-Sys}
\left\{\begin{array}{c}
\nabla \times \big(\mu^{-1} (\nabla \times E)  \big)- \omega^2 \eps  E =0 \mbox{ in } \Omega, \\[6pt]
\nabla \times \big(\hmu^{-1} (\nabla \times \hE)  \big)- \omega^2 \heps  \hE =0 \mbox{ in } \Omega, \\[6pt]
\hE \times \nu = E \times \nu  \mbox{ on } \partial \Omega, \\[6pt]
\big( \hmu^{-1 }(\nabla \times \hE) \big) \times \nu = \big( \mu^{-1 }(\nabla \times E) \big) \times \nu   \mbox{ on } \partial \Omega. 
\end{array} \right. 
\end{equation}
Then a complex number $\omega \in \C$ is called a \textit{transmission eigenvalue} if there exists a non-zero solution $(E, \hE) \in [L^2(\Omega)]^{6}$ of  \eqref{rem-Form-Sys}.  \Cref{thm1} might be translated as follows: 

{\it Completeness:} Assume that  $\eps, \, \mu, \, \heps, \, \hat{\mu}  \in [C^2(\bar \Omega)]^{3 \times 3}$ and \eqref{H} holds.   The space spanned by the generalized eigenfunctions  is complete in ${\bf G} (\Omega)$, i.e., the space spanned by them is dense in  ${\bf G}(\Omega)$, where 
\begin{multline}
{\bf G}(\Omega) = \Big\{ (u, \hu) \in [H(\curl, \Omega)]^2; \dive (\eps u) =\dive (\heps \hu) = 0 \mbox{ in } \Omega, \\[6pt]
 (\hu - u ) \times \nu = 0 \mbox{ on } \partial \Omega,  \, \big( \hmu^{-1}(\nabla \times \hu) \big) \times \nu - \big( \mu^{-1}(\nabla \times u) \big) \times \nu  = 0 \mbox{ on } \partial \Omega \Big\}
\end{multline}
\end{remark}

\begin{remark}  \rm 
In \cite{HM18}, the authors studied the completeness of generalized eigenfunctions in the isotropic case under the assumption that 
$$
\eps = \mu = \hmu = I \mbox{ in } \Omega, 
$$
$$
\heps \in C^\infty(\bar \Omega) \mbox{ and } \heps \mbox{ is constant different from 1 in a neighborhood of $\partial \Omega$}.  
$$
They considered the system under the form \eqref{rem-Form-Sys}.  Since $\eps = \mu = I$, their settings and ours are disjoint. Nevertheless, it seems that the boundary conditions  given in the definition of ${\bf G}(\Omega)$ are missing in \cite{HM18}. 
\end{remark}

\subsection{Proof of \Cref{thm1}} Applying  \Cref{pro-HS},  one has 
\begin{enumerate}
\item[-] $\cT_k^2: \bH(\Omega) \to \bH(\Omega)$ is a Hilbert-Schmidt operator.  
\item[-] for $\theta \in \mR$ with $|\Im (e^{2 i \theta})| > 0$, $e^{i \theta}$  is a direction of minimal growth of the modified resolvent of $\cT_k^2$.  
\end{enumerate}
Applying the theory of Hilbert-Schmidt operators, see e.g.  \cite[Theorem 16.4]{Agmon}, one derives that 

1) the closure of the space spanned by all generalized eigenfunctions of 
$\cT_{k}^{2}$ is equal to $\overline{\cT_k^2(\bH(\Omega))}$ (the closures are taken with respect to the $[L^2(\Omega)]^{12}$-norm).

\medskip
\noindent On the other hand, we have 

2) $\overline{\cT_k^2(\bH(\Omega))}=\bH(\Omega)$ since 
\[
\begin{aligned}
\bH(\Omega) &= \overline{\cT_k(\bH(\Omega))} && \mbox{(by Proposition \ref{pro-density})} \\
&= \overline{\cT_k \overline{\cT_k(\bH(\Omega))}} && \mbox{(by Proposition \ref{pro-density})} \\
&= \overline{\cT_k^2(\bH(\Omega))} && \mbox{(by the continuity of }\cT_k).
\end{aligned}
\]

3) The space spanned by the generalized eigenfunctions of  $\cT_{k}^{2}$ associated to the non-zero eigenvalues of $\cT_{k}^{2}$  is equal to the space spanned by the generalized  eigenfunctions of $\cT_k$ associated to the non-zero eigenvalues of $\cT_k$. This can be done as in the last part of the proof of \cite[Theorem 16.5]{Agmon}. Consequently,  the space spanned by all  generalized eigenfunctions of  $\cT_{k}^{2}$ is equal to  the space spanned by all  generalized  eigenfunctions of $\cT_k$. 

\medskip 
The conclusion now follows from 1), 2), and 3). \qed

\section{An upper bound for the counting function - Proof of \Cref{thm2}} \label{sect-thm2}

Let $\tlambda_j$ be the non-zero eigenvalues of $\cT_k$. Note that the non-zero eigenvalue values of $\cT_k^2$, counted according to multiplicity,  are $\tlambda_j^2$ (this can be proved as in the last part of the proof of \cite[Theorem 16.5]{Agmon}).  Applying the spectral theory of Hilbert-Schmidt operators, see e.g. \cite[Theorem 12.14]{Agmon} to $\cT_k^2$, we have 
\begin{equation}
\sum_{j} |\tlambda_j|^{4} \le \vvvert \cT_k^2 \vvvert^2. 
\end{equation}
Applying $i)$ of \Cref{pro-HS}, we obtain 
\[
\sum_{j}  |\tlambda_j|^{4} \leq C  |k|^{-1}.
\]
Note that  $\lambda_j$ is an transmission eigenvalue if and only if $(i \lambda_j - k)^{-1}$ is an eigenvalue of $\cT_k$, and they have the same multiplicity. It follows that 
\begin{equation}\label{thm2-p1}
\sum_{j}  \frac{1}{ |i \lambda_j - k|^{4}} \leq C  |k|^{-1}.
\end{equation}
Note that if $|\lambda_j| \leq |k|$, 
then  $|i \lambda_j - k| \leq 2|k|$. We then derive from \eqref{thm2-p1} that 
\[
\frac{1}{|k|^4}\sum_{j: |\lambda_j|\leq |k|} 1\leq C  |k|^{-1}.
\]
This implies 
\[
\cN(|k|) \le  C |k|^3.
\]
The proof is complete. \qed 

\providecommand{\bysame}{\leavevmode\hbox to3em{\hrulefill}\thinspace}
\providecommand{\MR}{\relax\ifhmode\unskip\space\fi MR }
\providecommand{\MRhref}[2]{%
  \href{http://www.ams.org/mathscinet-getitem?mr=#1}{#2}
}
\providecommand{\href}[2]{#2}


\begin{thebibliography}{10}

\bibitem{Agmon}
Shmuel Agmon, \emph{Lectures on elliptic boundary value problems}, Prepared for
  publication by B. Frank Jones, Jr. with the assistance of George W. Batten,
  Jr. Van Nostrand Mathematical Studies, No. 2, D. Van Nostrand Co., Inc.,
  Princeton, N.J.-Toronto-London, 1965. \MR{0178246}

\bibitem{ADNI}
Shmuel Agmon, Avron Douglis, and Louis Nirenberg, \emph{Estimates near the
  boundary for solutions of elliptic partial differential equations satisfying
  general boundary conditions. {I}}, Comm. Pure Appl. Math. \textbf{12} (1959),
  623--727. \MR{125307}

\bibitem{ADNII}
Shmuel Agmon, Avron Douglis, and Louis Nirenberg, \emph{{Estimates near the boundary for solutions of elliptic partial
  differential equations satisfying general boundary conditions. II.}}, Comm.
  Pure Appl. Math. \textbf{17} (1964), 35--92.

\bibitem{BCH11}
Anne-Sophie Bonnet-Ben~Dhia, Lucas Chesnel, and Houssem Haddar, \emph{On the
  use of {$T$}-coercivity to study the interior transmission eigenvalue
  problem}, C. R. Math. Acad. Sci. Paris \textbf{349} (2011), no.~11-12,
  647--651. \MR{2817384}

\bibitem{Brezis-FA}
Haim Brezis, \emph{Functional analysis, {S}obolev spaces and partial
  differential equations}, Universitext, Springer, New York, 2011. \MR{2759829}

\bibitem{CCH16}
Fioralba Cakoni, David Colton, and Houssem Haddar, \emph{Inverse scattering
  theory and transmission eigenvalues}, CBMS-NSF Regional Conference Series in
  Applied Mathematics, vol.~88, Society for Industrial and Applied Mathematics
  (SIAM), Philadelphia, PA, 2016. \MR{3601119}

\bibitem{CGH10}
Fioralba Cakoni, Drossos Gintides, and Houssem Haddar, \emph{The existence of
  an infinite discrete set of transmission eigenvalues}, SIAM J. Math. Anal.
  \textbf{42} (2010), no.~1, 237--255. \MR{2596553}

\bibitem{CHM15}
Fioralba Cakoni, Houssem Haddar, and Shixu Meng, \emph{Boundary integral
  equations for the transmission eigenvalue problem for {M}axwell's equations},
  J. Integral Equations Appl. \textbf{27} (2015), no.~3, 375--406. \MR{3435806}

\bibitem{Cakoni-Ng21}
Fioralba Cakoni and Hoai-Minh Nguyen, \emph{On the {D}iscreteness of
  {T}ransmission {E}igenvalues for the {M}axwell {E}quations}, SIAM J. Math.
  Anal. \textbf{53} (2021), no.~1, 888--913. \MR{4209665}

\bibitem{Chesnel12}
Lucas Chesnel, \emph{Interior transmission eigenvalue problem for {M}axwell's
  equations: the {$T$}-coercivity as an alternative approach}, Inverse Problems
  \textbf{28} (2012), no.~6, 065005, 14. \MR{2924302}

\bibitem{CM88}
David Colton and Peter Monk, \emph{The inverse scattering problem for
  time-harmonic acoustic waves in an inhomogeneous medium}, Quart. J. Mech.
  Appl. Math. \textbf{41} (1988), no.~1, 97--125. \MR{934695}

\bibitem{Gagliardo59}
Emilio Gagliardo, \emph{Ulteriori propriet\`a di alcune classi di funzioni in
  pi\`u variabili}, Ricerche Mat. \textbf{8} (1959), 24--51. \MR{109295}

\bibitem{GR86}
Vivette Girault and Pierre-Arnaud Raviart, \emph{Finite element methods for
  {N}avier-{S}tokes equations}, Springer Series in Computational Mathematics,
  vol.~5, Springer-Verlag, Berlin, 1986, Theory and algorithms. \MR{851383}

\bibitem{Haddar04}
Houssem Haddar, \emph{The interior transmission problem for anisotropic
  {M}axwell's equations and its applications to the inverse problem}, Math.
  Methods Appl. Sci. \textbf{27} (2004), no.~18, 2111--2129. \MR{2102315}

\bibitem{HM18}
Houssem Haddar and Shixu Meng, \emph{The spectral analysis of the interior
  transmission eigenvalue problem for {M}axwell's equations}, J. Math. Pures
  Appl. (9) \textbf{120} (2018), 1--32. \MR{3906154}

\bibitem{KIrsch86}
Andreas Kirsch, \emph{The denseness of the far field patterns for the
  transmission problem}, IMA J. Appl. Math. \textbf{37} (1986), no.~3,
  213--225. \MR{983987}

\bibitem{LV12}
Evegeny Lakshtanov and Boris Vainberg, \emph{Ellipticity in the interior
  transmission problem in anisotropic media}, SIAM J. Math. Anal. \textbf{44}
  (2012), no.~2, 1165--1174. \MR{2914264}

\bibitem{Ng-Complementary}
Hoai-Minh Nguyen, \emph{Asymptotic behavior of solutions to the {H}elmholtz
  equations with sign changing coefficients}, Trans. Amer. Math. Soc.
  \textbf{367} (2015), no.~9, 6581--6595. \MR{3356948}

\bibitem{Ng-WP}
Hoai-Minh Nguyen,  \emph{Limiting absorption principle and well-posedness for the
  {H}elmholtz equation with sign changing coefficients}, J. Math. Pures Appl.
  (9) \textbf{106} (2016), no.~2, 342--374. \MR{3515306}

\bibitem{Ng-Superlensing-Maxwell}
Hoai-Minh Nguyen,  \emph{Superlensing using complementary media and reflecting
  complementary media for electromagnetic waves}, Adv. Nonlinear Anal.
  \textbf{7} (2018), no.~4, 449--467. \MR{3871415}

\bibitem{Ng-CALR-O-M}
Hoai-Minh Nguyen,  \emph{Cloaking property of a plasmonic structure in doubly
  complementary media and three-sphere inequalities with partial data},
  (2019), preprint, https://arxiv.org/abs/1912.09098.

\bibitem{Ng-Negative-Cloaking-M}
Hoai-Minh Nguyen,  \emph{Cloaking using complementary media for electromagnetic waves},
  ESAIM Control Optim. Calc. Var. \textbf{25} (2019), Art. 29, 19. \MR{3990650}

\bibitem{Ng-CALR-M}
Hoai-Minh Nguyen,  \emph{The invisibility via anomalous localized resonance of a source
  for electromagnetic waves}, Res. Math. Sci. \textbf{6} (2019), no.~4, Paper
  No. 32, 22. \MR{4011564}

\bibitem{MinhHung17}
Hoai-Minh Nguyen and Quoc-Hung Nguyen, \emph{Discreteness of interior
  transmission eigenvalues revisited}, Calc. Var. Partial Differential
  Equations \textbf{56} (2017), no.~2, Paper No. 51, 38. \MR{3626617}

\bibitem{MinhHung2}
Hoai-Minh Nguyen and Quoc-Hung Nguyen, \emph{The {W}eyl law of transmission eigenvalues and the completeness
  of generalized transmission eigenfunctions},  (2020),
  https://arxiv.org/abs/2008.08540.

\bibitem{NgSil}
Hoai-Minh Nguyen and Swarnendu Sil, \emph{Limiting {A}bsorption {P}rinciple and
  {W}ell-{P}osedness for the {T}ime-{H}armonic {M}axwell {E}quations with
  {A}nisotropic {S}ign-{C}hanging {C}oefficients}, Comm. Math. Phys.
  \textbf{379} (2020), no.~1, 145--176. \MR{4152269}

\bibitem{Nirenberg59}
Louis Nirenberg, \emph{On elliptic partial differential equations}, Ann. Scuola
  Norm. Sup. Pisa Cl. Sci. (3) \textbf{13} (1959), 115--162. \MR{109940}

\bibitem{Robbiano13}
Luc Robbiano, \emph{Spectral analysis of the interior transmission eigenvalue
  problem}, Inverse Problems \textbf{29} (2013), no.~10, 104001, 28.
  \MR{3116196}

\bibitem{Robbiano16}
Luc Robbiano, \emph{Counting function for interior transmission eigenvalues}, Math.
  Control Relat. Fields \textbf{6} (2016), no.~1, 167--183. \MR{3448675}

\bibitem{BF07}
Mikhail Shl\"emovich and Nikolai~D. Filonov, \emph{Weyl asymptotics of the
  spectrum of the {M}axwell operator with non-smooth coefficients in
  {L}ipschitz domains}, Nonlinear equations and spectral theory, Amer. Math.
  Soc. Transl. Ser. 2, vol. 220, Amer. Math. Soc., Providence, RI, 2007,
  pp.~27--44. \MR{2343605}

\bibitem{Sylvester12}
John Sylvester, \emph{Discreteness of transmission eigenvalues via upper
  triangular compact operators}, SIAM J. Math. Anal. \textbf{44} (2012), no.~1,
  341--354. \MR{2888291}

\bibitem{Vodev15}
Georgi Vodev, \emph{Transmission eigenvalue-free regions}, Comm. Math. Phys.
  \textbf{336} (2015), no.~3, 1141--1166. \MR{3324140}

\bibitem{Vodev18}
Georgi Vodev, \emph{High-frequency approximation of the interior
  {D}irichlet-to-{N}eumann map and applications to the transmission
  eigenvalues}, Anal. PDE \textbf{11} (2018), no.~1, 213--236. \MR{3707296}

\bibitem{Vodev21}
Georgi Vodev, \emph{Semiclassical parametrix for the {M}axwell equation and
  applications to the electromagnetic transmission eigenvalues},  (2021),
  https://arxiv.org/abs/2102.08662.

\bibitem{Weyl12}
H.~Weyl, \emph{\"{U}ber das {S}pektrum der {H}ohlraumstrahlung}, J. Reine
  Angew. Math. \textbf{141} (1912), 163--181. \MR{1580849}

\end{thebibliography}
\end{document}